\documentclass[12pt]{article}
\usepackage[final]{epsfig}
\usepackage{graphics}
\usepackage{amsmath}
\usepackage{amsfonts}
\usepackage{latexsym}
\usepackage{amssymb}
\usepackage{graphicx}
\usepackage{epstopdf}

\newtheorem{lemma}{Lemma}[section]
\newtheorem{proposition}[lemma]{Proposition}
\newtheorem{remark}[lemma]{Remark}
\newtheorem{example}[lemma]{Example}
\newtheorem{theorem}{Theorem}

\newtheorem{conjecture}[theorem]{Conjecture}

\begin{document}
\newcommand{\eps}{{\varepsilon}}
\newcommand{\proofend}{$\Box$\bigskip}
\newcommand{\C}{{\mathbf C}}
\newcommand{\Q}{{\mathbf Q}}
\newcommand{\R}{{\mathbf R}}
\newcommand{\Z}{{\mathbf Z}}
\newcommand{\RP}{{\mathbf {RP}}}

\newcommand\pd[2]{\frac{\partial #1}{\partial #2}}
\def\proof{\paragraph{Proof.}}

\title {On configuration spaces of plane polygons, sub-Riemannian geometry and periodic orbits of outer billiards}
\author{Daniel Genin and Serge Tabachnikov\\
{\it Department of Mathematics, Penn State University}\\
{\it University Park, PA 16802, USA}
}
\date{}
\maketitle

\section{Introduction} \label{intro}

The classical billiard system describes the motion of a point in a plane domain subject to the elastic reflection off the boundary, described by the familiar law of geometrical optics: the angle of incidence equals the angle of reflection; see, e.g.,  \cite{Tab1,Tab2} for surveys of mathematical billiards.

For every $n\geq 2$, the billiard system inside a circle has a very special property:  every point of the circle is the starting point of an $n$-periodic billiard orbit; this orbit is an inscribed regular $n$-gon.\footnote{Similarly, one can consider periodic orbits with other rotation numbers corresponding to regular star-shaped polygons.} Likewise, an ellipse enjoys the same property for every $n\geq 3$; the orbits are inscribed $n$-gons of extremal perimeter length.\footnote{For circles and ellipses, this behavior is due to the complete integrability of the billiard ball map.}
A billiard table of constant width also has the property that every point on its boundary belongs to a $2$-periodic, back-and-forth, billiard trajectory, and this is a dynamic characterization of the  curves of constant width. The phase space of the billiard ball map is a cylinder, and the family of $n$-periodic orbits forms an invariant circle of the billiard ball map consisting of $n$-periodic points. 

How exceptional is this property? More specifically, given $n\geq 3$, one wants to describe plane billiards such that the billiard ball map has an invariant curve consisting of $n$-periodic points. The first result in this direction was obtained by Innami \cite{Inn}: he constructed a billiard curve, other than a circle or an ellipse, whose every point is  a vertex of a $3$-periodic billiard trajectory (these billiard triangles are isosceles in his example).\footnote{A similar result was obtained by M. Berger, unpublished.}

In  a recent paper \cite{BZh}, Baryshnikov and Zharnitskii developed a new elegant approach to this problem. In a nutshell, their innovation was in switching focus from billiard curves to the space of plane $n$-gons. If such an $n$-gon, $z_1,\dots,z_n$, is a periodic billiard trajectory then, due to the law ``the angle of incidence equals the angle of reflection", we know the directions of the billiard curve at each vertex $z_i$. This determines an $n$-dimensional distribution ${\cal D}$ in the space of $n$-gons, called the {\it Birkhoff distribution}. One can further restrict attention to the hypersurface  of $n$-gons with a fixed perimeter length: the Birkhoff distribution is everywhere tangent to these hypersurfaces.

It turns out that the Birkhoff distribution is non-integrable. A 1-parameter family of billiard $n$-gons $z_1(t),\dots,z_n(t),\ 0\leq t\leq 1$, is interpreted as a curve in the space of $n$-gons, tangent to ${\cal D}$ and satisfying the ``monodromy" condition $z_i(1)=z_{i+1}(0),\ i=1,\dots,n$. The problem of finding such {\it horizontal} curves belongs to sub-Riemannian geometry and non-holonomic dynamics. 

The main result of \cite{BZh} is that, for every $n\geq 2$, the space of billiard tables, for which the billiard ball map possesses an invariant curve consisting of $n$-periodic points, is infinite-dimensional (this space also has an infinite codimension in the space of all billiard tables).

\begin{figure}[hbtp]
\centering
\includegraphics[width=2in]{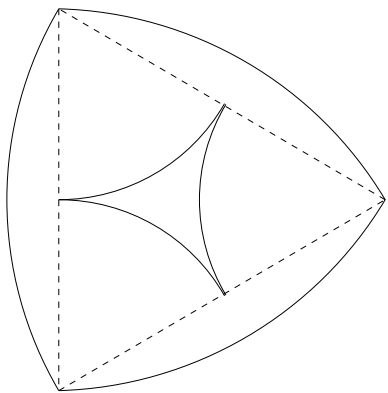}
\caption{Reuleaux triangle and the curve $\Delta$}
\label{delta}
\end{figure}

\begin{example}
{\rm
Let us illustrate this approach in the simplest case  $n=2$, that is, curves of constant width (for information on billiards of constant width, see \cite{Kn}). 

Consider the space of  oriented segments, say, of length 2.  A segment $z_1 z_2$ is ``allowed" to move in such a way that the velocities of both end-points are orthogonal to the segment. Let $O=(x,y)$ be the mid-point of the segment and $\alpha$ its direction. Then the space of segments has coordinates $(x,y,\alpha)$ and
$$
z_1=(x-\cos \alpha, y-\sin \alpha),\ z_2=(x+\cos \alpha, y+\sin \alpha).
$$
We conclude that the velocity of point $O$ must be orthogonal to the vector $(\cos \alpha, \sin \alpha)$, and this is the definition of the Birkhoff distribution in this case. Thus ${\cal D}$ is the kernel of the 1-form $\lambda=\cos \alpha\ dx + \sin \alpha\ dy$. This 1-form is contact: $\lambda \wedge d\lambda \neq 0$. We have identified our space of segments with the space of cooriented contact elements in the plane, a fundamental example of a contact 3-dimensional manifold (see, e.g., \cite{Ge}).

If $z_1(t) z_2(t)$ is a curve, tangent to ${\cal D}$ and satisfying the monodromy condition $z_1(1)=z_2(0), z_2(1)=z_1(0)$, then the respective curve $(x(t),y(t),\alpha(t))$ is Legendrian (i.e., tangent to the contact distribution) and satisfies $x(1)=x(0), y(1)=y(0)$ and $\alpha(1)=\alpha(0)+\pi$. 

The projection of this Legendrian curve on the $(x,y)$-plane, that is, the trajectory of the mid-point $O$, is a closed smooth curve $\Delta$ with singularities (generically, semi-cubical cusps); the total winding of this curve is $\pi$.
Every such curve $\Delta$ uniquely lifts to a Legendrian curve in the space of contact elements: the missing $\alpha$ coordinate is recovered  as the slope: $\cot \alpha=-dy/dx$. 

In other words, we use $\Delta$ as a guide: place a segment of length 2 so that its mid-point $O$ is on the curve $\Delta$ and the segment is orthogonal to $\Delta$, and then slide  the segment  around $\Delta$. 
If certain convexity conditions hold (for example, $\Delta$ should not have inflection points) then the  end-points $z_1$ and $z_2$ will describe a curve of constant width.

For example,  if the end-points $z_1$ and $z_2$ move along a circle then $\Delta$ degenerates to its center, and if $z_1$ and $z_2$ describe
the Reuleaux triangle then $\Delta$ is  made of three $60^{\circ}$ arcs of a circle, see figure \ref{delta}.
}
\end{example}

In this paper we extend the results from \cite{BZh} to outer (or dual) billiards, a dynamical system that is a close relative of the conventional, inner, billiards; see \cite{D-T,Tab1,Tab2} for surveys and \cite{Gen} for a recent interesting result.

\begin{figure}[hbtp]
\centering
\includegraphics[width=2.5in]{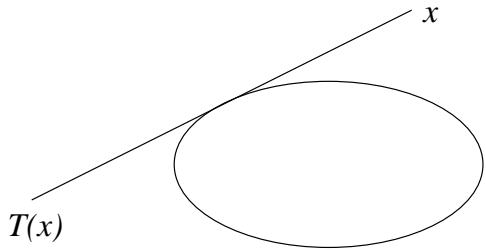}
\caption{Outer billiard map}
\label{outb}
\end{figure}

An outer billiard table  is a strictly convex compact plane domain $D$.  Pick a point $x$ outside $D$. There are two support lines from $x$ to $D$; choose  one of them, say, the right one from $x$'s view-point, and reflect $x$ in the support point. One obtains a new point, $y$, and the transformation $T: x\mapsto y$ is the outer billiard map, see figure \ref{outb}. If the boundary of the outer billiard table contains straight segments (for example, is a convex polygon), the map $T$ or its inverse is not defined on the lines containing these segments; this still leaves a set of full measure on which all iterations of $T$ are defined. We will consider only outer billiard tables with strictly convex smooth boundaries.

Let us mention four fundamental properties of outer billiards, important for us. First, 
the outer billiard map commutes with affine transformations of the plane. Secondly,
the outer billiard map preserves the standard area form in the plane (similarly to the billiard ball map which is an area preserving map of a cylinder). Thirdly, the outer billiard map is a twist map of the exterior of the outer billiard table $D$, a semi-infinite cylinder, with respect to the vertical foliation consisting of tangent half-lines to the boundary of $D$. And finally, an $n$-periodic outer billiard  orbit is an  $n$-gon (possibly, star-shaped), circumscribed about the table $D$ and having an extremal area. This is an outer analog of a variational property  of periodic billiard trajectories: they are inscribed polygons of extremal perimeter length. As a consequence, if an outer billiard has a 1-parameter family of $n$-periodic orbits then the respective circumscribed  $n$-gons have the same areas.

Let $G_n$ be the {\it cyclic configuration space} of the plane, that is, the space of polygons $z_1,\dots, z_n$ such that $z_i \neq z_{i+1}$ for $i=1,\dots,n$ (where the index $i$ is understood cyclically). 
If $z_1(t),\dots, z_n(t)$ is a 1-parameter family of $n$-periodic outer billiard orbits with $T(z_i(t))=z_{i+1}(t)$  then each segment $z_i(t) z_{i+1}(t)$ envelopes the boundary of the outer billiard table and is bisected by the tangency point. This prompts the following definition.

Denote by ${\cal F}$ the $n$-dimensional distribution in $G_n$ determined by the following condition: the velocities of points $z_1,\dots, z_n$ lie in ${\cal F}$ if,
for each $i=1,\dots,n$, the motion of the line $z_i z_{i+1}$ is an infinitesimal rotation about the mid-point of the segment $z_i z_{i+1}$. We call ${\cal F}$ the {\it dual Birkhoff distribution}.

The  dual Birkhoff distribution is tangent to the level hypersurfaces of the area function, and we can restrict attention to non-degenerate polygons of a fixed non-zero area. On this hypersurface, ${\cal F}$ is non-integrable: the vector fields, tangent to ${\cal F}$, and their first commutators generate the tangent space at every point (see Theorem \ref{nonint}). 

We are interested in the outer billiard tables possessing invariant curves consisting of $n$-periodic orbits. Such an invariant curve gives rise to a horizontal curve $z_1(t),\dots, z_n(t),\ 0\leq t\leq 1 $, in $G_n$ satisfying the monodromy condition $z_i(1)=z_{i+1}(0),\ i=1,\dots, n$. If the outer billiard table is a circle (or an ellipse) then, for every  pair $(n,k)$ where  $3\leq n, 1\leq k\leq n/2$, and $k$ coprime with $n$, there is an invariant curve consisting  of $n$-periodic orbits with rotation number $k$ (that is, making $k$ turns around $D$). We show that, for every such pair $(n,k)$, one can perturb a circle so that the resulting outer billiard still has an invariant curve consisting of $n$-periodic orbits with rotation number $k$; the space of such deformations is  infinite-dimensional and depends on functional parameters (see Theorem \ref{mfld}).

In Section \ref{per3}, we study the simplest case of $n=3$ in detail; this is the closest outer billiard analog of curves of constant width. This case illustrates the general approach by explicit computations. We also discuss a discretization of the problem in which a smooth outer billiard curve is replaced by a polygon.

We also prove that, for every outer billiard table, the set of $3$-periodic orbits has zero measure (see Theorem \ref{3per}). For inner billiards, this is a well-known result having a number of different proofs \cite{Ry,Sto,Vo,Wo}. We also give different proofs, one of them follows ideas of \cite{Wo}  and another those of \cite{BZh}. Note that there exist
area preserving twist maps of a cylinder whose set of periodic points has a non-empty interior.

Let us finish this introduction with a problem: {\it can invariant billiard curves (inner or outer) consisting of periodic points with different periods, or the same period but different rotation numbers, coexist, unless the table is an ellipse? In particular, is it true that 
the only billiard of constant width, having an invariant curve consisting of 3-periodic points, is a circle?}

\section{Configuration space of plane polygons and the dual Birkhoff distribution} \label{conf}

Let us fix some notations. The determinant of two vectors, that is, the cross-product, is denoted by $[\ ,\ ]$. A generic point in the plane is denoted by $z=(x,y)$. If $z_1=(x_1,y_1)$ and $z_2=(x_2,y_2)$ then $z_1 dz_2$  denotes the 1-form $x_1 dx_2+y_1 dy_2$. Likewise, $[z_1,dz_2]$  denotes the 1-form $x_1 dy_2-y_1 dx_2$ and $[dz,dz]$ the 2-form $2dx \wedge dy$.
If $w=(u,v)$ is a vector in $\R^2$ and $z=(x,y)$ is a point then by $w \partial z$ we mean the tangent vector $u \partial x + v \partial y$.  When working with $n$-gons, we always understand the indices cyclically, so that $i=i+n$. For a polygon $z_1,\dots,z_n$, set
$$
A(z_1,\dots,z_n)=\sum_{i=1}^n [z_i,z_{i+1}];
$$
this is twice the  area of the polygon (counted with appropriate multiplicities and signs).

Let ${\cal F}$ be the dual Birkhoff distribution on the space of $n$-gons defined in Section \ref{intro}. Let $Z=(z_1,\dots,z_n) \in G_n$ be an $n$-gon and $W=(w_1,\dots,w_n)$ a tangent vector to $G_n$ at $Z$; the vector $w_i$ is the velocity of the vertex $z_i$.

\begin{lemma} \label{chardist}
${\cal F}$ is an $n$-dimensional distribution. 
A vector $W$ lies in ${\cal F}$ if and only if, for all $i=1,\dots,n$, one has: $[w_i+w_{i+1},z_{i+1}-z_i]=0$, that is, the vectors $w_i+w_{i+1}$ and $z_{i+1}-z_i$ are collinear.
\end{lemma}

\proof To see that ${\cal F}$ is  $n$-dimensional, view a polygon 
$Z=(z_1,\dots,z_n) \in G_n$ as an ordered collection of $n$ lines $L=(z_1z_2,z_2z_3,\dots,z_nz_1)$. From this dual viewpoint, $G_n$ becomes an open subset in the space ${\cal L}_n$ of $n$-tuples of oriented  lines in the plane. An infinitesimal rotation of the $i$-th line $z_i z_{i+1}$ about the mid-point of the segment $z_i z_{i+1}$ determines a tangent vector to ${\cal L}_n$ at point $L$, and these vectors are clearly independent for $i=1,\dots,n$. It follows that dim ${\cal F}=n$.

For $W$ to belong to ${\cal F}$, the line through points $z_i+\eps w_i$ and $z_{i+1}+\eps w_{i+1}$ must pass through the mid-point $(z_i + z_{i+1})/2$; here $\eps$ is infinitesimal. Thus the linear term in $\eps$ of the determinant 
$$
\left [\frac{z_{i+1}-z_i}{2} +\eps w_{i+1}, \frac{z_{i}-z_{i+1}}{2} +\eps w_{i}\right]
$$ vanishes, which is equivalent to $[w_i+w_{i+1},z_{i+1}-z_i]=0$.
\proofend

An $n$-gon is called {\it non-degenerate} if no three consecutive vertices  lie on one line. 
The set of non-degenerate $n$-gons is denoted by $U_n$. Clearly, $U_n$ is an open subset of $G_n$.

Let us introduce the 1-forms $\alpha_i=[z_{i+1}-z_i, dz_{i+1}+dz_i],\ i=1,\dots,n$. 

\begin{lemma} \label{forms}
Let $Z\in U_n$ be a non-degenerate polygon.  Then the 1-forms $\alpha_i$ are linearly independent at $Z$ and the intersection  of their kernels is the fiber $F(Z)$ of
${\cal F}$ at point $Z$.
\end{lemma}

\proof First, we claim that each form $\alpha_i$ vanishes on every tangent vector $W \in {\cal F}$. Indeed, $\alpha_i (W)=[z_{i+1}-z_i,w_i+w_{i+1}]$, and this is zero, according to Lemma \ref{chardist}.

Suppose that a non-trivial linear combination $\sum t_i \alpha_i$ vanishes at point $Z$. Then, for every test vector $W=(w_1,\dots,w_n)$, one has:
$$
0=\sum t_i [z_{i+1}-z_i,w_i+w_{i+1}] = \sum t_i  [z_{i+1}-z_i,w_i] + \sum t_i [z_{i+1}-z_i,w_{i+1}] = 
$$
$$
\sum[t_i (z_{i+1}-z_i),w_i] + \sum [t_{i-1}(z_i-z_{i-1}),w_i]= \sum[t_i (z_{i+1}-z_i)+t_{i-1}(z_i-z_{i-1}),w_i].
$$
Therefore  
\begin{equation} \label{lindep}
t_i (z_{i+1}-z_i)+t_{i-1}(z_i-z_{i-1})=0,\ \ i=1,\dots,n.
\end{equation}
 Since the polygon $Z$ is non-degenerate,  the vectors  $z_{i+1}-z_i$ and $z_i-z_{i-1}$ are linearly independent, and (\ref{lindep}) implies that $t_i=0$ for all $i$.
 
Since  ${\cal F}$ is $n$-dimensional and the 1-forms $\alpha_i$ are linearly independent, one has: $\cap_{i=1}^n\ {\rm Ker}\ \alpha_i = {\cal F}$. 
\proofend

The next lemma shows that one can restrict attention to the level hypersurfaces of the area function.

\begin{lemma} \label{hypers}
Every non-zero level hypersurface of the function $A$ is smooth, and the distribution ${\cal F}$ is tangent to these hypersurfaces.
\end{lemma}

\proof Let $Q:V\to V^*$ be a self-adjoint linear operator in a vector space $V$. Consider the quadratic function $F(x)=Q(x)\cdot x$ on $V$. Then every non-zero level hypersurface of $F$ is smooth. Indeed, by Euler's formula, $\nabla F(x) \cdot x=2F(x)$, hence  $dF$ does not vanish on a non-zero level hypersurface.
The function $A$ is a quadratic form on the space of $n$-gons, and the first claim follows.

One has the  equality
$$
\sum_{i=1}^n \alpha_i = \sum_{i=1}^n [z_{i+1},dz_i]- [z_{i},dz_{i+1}]=-dA,
$$
which implies the second claim.
\proofend

\begin{figure}[hbtp]
\centering
\includegraphics[width=2.5in]{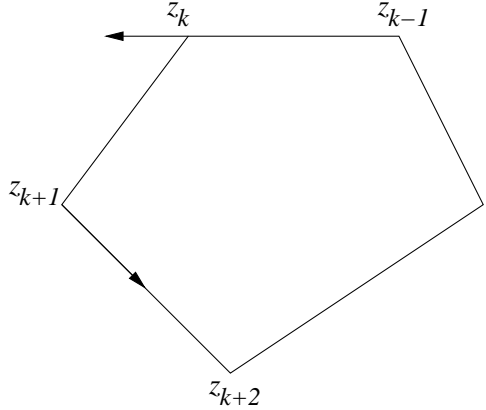}
\caption{Tangent vector $W_k$}
\label{vct}
\end{figure}

Let us construct a system of $n$ linearly independent vector fields tangent to ${\cal F}$. Let $Z=(z_1,\dots,z_n)$ be a non-degenerate polygon. Set $a_i=[z_i-z_{i-1},z_{i+1}-z_i]$,  twice the oriented area of the triangle $z_{i-1} z_i z_{i+1}$.
Denote by $W_k=(w_{1k},\dots,w_{nk})$ the tangent vector to $U_n$ at $Z$ whose $i$-th component is zero, unless $i=k$ or $i=k+1$. For these values of $i$, one has:
$$
w_{k,k}=a_{k+1} (z_k-z_{k-1}),\ \ w_{k+1,k}=a_{k} (z_{k+2}-z_{k+1}),
$$
see figure \ref{vct}.

\begin{lemma} \label{span}
The vectors $W_1,\dots,W_n$ are linearly independent at every point $Z \in U_n$ and they span $F(Z)$.
\end{lemma}

\proof
To show that $W_k$ lies in ${\cal F}$, it suffices to check that the vector $a_{k+1} (z_k-z_{k-1}) + a_{k} (z_{k+2}-z_{k+1})$ is collinear with the vector $z_{k+1}-z_k$. Indeed,
$$
[z_{k+1}-z_{k},z_{k+2}-z_{k+1}] (z_k-z_{k-1}) +  [z_k-z_{k-1},z_{k+1}-z_k] (z_{k+2}-z_{k+1})=
$$
$$
[z_k-z_{k-1},z_{k+2}-z_{k+1}] (z_{k+1}-z_{k}),
$$
which follows from the general identity, valid for any three vectors:
$$
[u,v]w+[v,w]u+[w,u]v=0.
$$
Under identification of the space of $n$-gons with the space of $n$-tuples of lines ${\cal L}_n$, the tangent vector $W_k$ corresponds to a nontrivial infinitesimal rotation of the line $z_k z_{k+1}$ about the mid-point of the segment  $z_k z_{k+1}$. Therefore these vectors are linearly independent and hence span the dual Birkhoff distribution. 
\proofend

Now we prove that the the dual Birkhoff distribution is completely non-integrable.

\begin{theorem} \label{nonint}
Let $M^{2n-1}$ be the hypersurface in $U_n$ consisting of non-degenerate $n$-gons with a fixed non-zero area. Then the tangent space of $M$ at every point is generated by the vectors fields tangent to ${\cal F}$ and their first commutators. Thus the bracket growth type of ${\cal F}$ is $(n,2n-1)$ (just like that of the Birkhoff distribution, see \cite{BZh}).
\end{theorem}

\proof 
Consider the vector fields $\xi_k=[W_{k-1}, W_k]$ (the fields $W_i$ and $W_j$ clearly  commute if $|i-j|\geq 2$). We want to show that these fields, along with ${\cal F}$, span the tangent space to $M$ at every point.

According to Lemmas \ref{forms} and \ref{hypers}, the 1-forms $\alpha_1,\dots,\alpha_{n-1}$ constitute a frame in the conormal space to ${\cal F}$ at every point $Z \in M$. For $k=1,\dots,n$, consider the vector $(\alpha_1(\xi_k),\dots,\alpha_{n-1}(\xi_k))\in \R^{n-1}$. It suffices to show that these $n$ vectors span $\R^{n-1}$.

Since the fields $W_k$ belong to the kernel of each form $\alpha_i$, the Cartan formula implies: $\alpha_i (\xi_k) = d\alpha_i (W_{k-1}, W_k)$. 

Recall that $Z=(z_1,\dots,z_n)$ with $z_i=(x_i,y_i)$. Let $\omega_i =dx_i \wedge dy_i,\ i=1,\dots,n$. We claim that $d\alpha_i = 2(\omega_{i+1}-\omega_i).$
Indeed, 
$$
d\alpha_i=[dz_{i+1}-dz_i,dz_{i+1}+dz_i]=[dz_{i+1},dz_{i+1}]-[dz_{i},dz_{i}]=2(\omega_{i+1}-\omega_i).
$$
Thus we want to show that the vectors 
$$
((\omega_{2}-\omega_1)(W_{k-1}, W_k),\dots,(\omega_{n}-\omega_{n-1})(W_{k-1}, W_k))
$$ 
span  $\R^{n-1}$ for $k=1,\dots,n$. This will follow if we show that the vectors
$$
\omega_1(W_{k-1}, W_k), \omega_2(W_{k-1}, W_k), \dots, \omega_n (W_{k-1}, W_k)
$$
span $\R^{n}$ for $k=1,\dots,n$. But $\omega_i(W_{k-1}, W_k)=0$, unless $i=k$, and 
$\omega_k (W_{k-1}, W_k)= a_{k-1}a_k a_{k+1} \neq 0$. Thus $\omega_i(W_{k-1}, W_k)$ is a non-degenerate diagonal matrix, and we are done.
\proofend

\begin{remark} \label{mid}
{\rm The dual Birkhoff distribution ${\cal F}$ on the cyclic configuration space $G_n$ is Lagrangian with respect to the following symplectic structure. The space of oriented lines in the plane, topologically, a cylinder, has a canonical area form (see, e.g., \cite{Tab1,Tab2}), and the direct product of these forms is a symplectic structure in the space of $n$-tuples of oriented lines ${\cal L}_n$. In the proof of Lemma \ref{chardist}, we realized $G_n$ as an open subset in ${\cal L}_n$, and this provides a symplectic structure on $G_n$. Infinitesimal rotations of distinct lines determine symplectically orthogonal tangent vectors to ${\cal L}_n$, hence ${\cal F}$ is Lagrangian. The Lagrangian property of the dual Birkhoff distribution echoes that of the Birkhoff distribution \cite{BZh}.
}
\end{remark}

\section{Outer billiards with invariant curves consisting of periodic points} \label{perpts}

In this section we show that, for every pair $(n,k)$ with $3\leq n, 1\leq k\leq n/2$ and $k$ coprime with $n$,  there is an abundance of outer billiard tables, possessing  invariant curves consisting of $n$-periodic orbits with rotation number $k$; such tables are constructed as perturbations of circles. 

Denote by $H^2_{\cal F}$ the space of parameterized horizontal paths $[0,1]\to M$ (paths tangent to the dual Birkhoff distribution) whose first 2 derivatives are square integrable (one can work with a higher number of derivatives as well). For a polygon $Z\in M$, denote by $H^2_{\cal F} (Z)$ the space of paths starting at $Z$. Then $H^2_{\cal F} (Z)$ is a Hilbert manifold, see \cite{Mon}. 

Let $Z$ be a regular star-shaped $n$-gon with rotation number $k$. 
We have the following result.

\begin{theorem} \label{mfld} 
The set of (parameterized) outer billiard curves, sufficiently close to a circular one and possessing  invariant curves consisting of  $n$-periodic orbits with rotation number $k$, is a smooth Hilbert submanifold of codimension $2n-1$ in $H^2_{\cal F} (Z)$.
\end{theorem}

\proof 
To fix ideas, consider  the case of simple periodic $n$-gons, that is, the case of $k=1$. At the end of the proof, we will indicate the small adjustment  to be made in the case of $k>1$.

Start with a circular outer billiard. One has a 1-parameter family of periodic $n$-gons 
\begin{equation} \label{1curve}
Z(t)=(z_1(t),\dots,z_n(t)) \ \ {\rm with}\  \ z_j(t)=e^{2\pi i(j+t)/n},\ t\in[0,1]
\end{equation}
circumscribed about the circular table (the table is scaled appropriately). Let $\sigma: G_n \to G_n$ be the cyclic permutation of the vertices of a polygon: $\sigma(z_1,\dots,z_n)=(z_2,\dots,z_n,z_1).$ Then $Z(1)=\sigma(Z(0))$.

The family (\ref{1curve}) determines a horizontal curve $\gamma_0(t)$ in $G_n$ whose end-points are $Z(0)$ and $Z(1)$. 
Denote by $H^2_{\cal F} (Z(0),Z(1))\subset H^2_{\cal F} (Z(0))$ the set of horizontal curves with the terminal point $Z(1)$.  We want to perturb $\gamma_0(t)$ in $H^2_{\cal F} (Z(0),Z(1))$. If a perturbation is sufficiently small and trivial in a neighborhood of end-points then we obtain a 1-parameter family of convex $n$-gons whose vertices trace $n$ arcs that smoothly match together to form an invariant curve of the outer billiard map consisting of $n$-periodic points; the outer billiard curve is recovered as the envelope of the segments connecting the consecutive vertices. 

We want to show that $H^2_{\cal F} (Z(0),Z(1))$ is a  Hilbert submanifold of codimension $2n-1$ in $H^2_{\cal F} (Z(0))$. Denote by $E: H^2_{\cal F} (Z(0)) \to M$ the end-point map that assigns the end-point $\gamma(1)$ to a path $\gamma$. This map is smooth, and if  $E$ is a submersion at  $\gamma_0 \in H^2_{\cal F} (Z(0))$ then $H^2_{\cal F} (Z(0),Z(1))=E^{-1}(Z(1))$ is a codimension $2n-1$ submanifold, see \cite{Mon}. A curve which is a singular point of the end-point map is called {\it singular}.

To prove that $\gamma_0$ is not singular we use a criterion due to Hsu \cite{Hsu}. Denote by ${\cal F}^{\perp}\subset T^* M$ the annihilator of the dual Birkhoff distribution and let $\Omega$ be the restriction of the canonical symplectic form on the cotangent bundle to ${\cal F}^{\perp}$. Let $\pi:T^* M\to M$ be the projection. A characteristic in ${\cal F}^{\perp}$ is a curve that does not intersect the zero section and whose direction at every point lies in the kernel of $\Omega$. By Hsu's criterion, a horizontal curve is singular if and only it is the projection under $\pi$ of a characteristic curve in ${\cal F}^{\perp}$. 

It is convenient for computations to use the following form of this criterion. The 1-forms $\alpha_i$, introduced in Section \ref{conf},  generate ${\cal F}^{\perp}$: every covector in ${\cal F}^{\perp}$ can be written as $\sum \lambda_i \alpha_i$ where the coefficients $\lambda_i$ are ``Lagrange multipliers". There is exactly one relation: $\sum \alpha_i =0$  (see Lemma \ref{hypers}), and to factor this relation out, we assume that $\sum \lambda_i =0$. Consider the 2-form $\omega(\lambda)=\sum \lambda_i d\alpha_i$ on $M$. Then one has the following result, see \cite{Mon}: the kernel of $\Omega$ at a point  $(x,\sum \lambda_i \alpha_i) \in {\cal F}^{\perp}$ is isomorphically projected by the differential $d\pi$ onto the kernel of the form $\omega(\lambda)$ restricted to the fiber $F(x)$ at point $x \in M$.

Let us find the 2-form $\omega(\lambda)$. From the proof of Theorem \ref{nonint}, we know that $d\alpha_i=2(\omega_{i+1}-\omega_i)$ where $\omega_i=dx_i \wedge dy_i$. Thus 
\begin{equation} \label{ker1}
-\frac{1}{2}\omega(\lambda)=\sum_i \lambda_i (\omega_i -\omega_{i+1})=\sum_i(\lambda_i - \lambda_{i-1}) \omega_i.
\end{equation}
Assume that $\dot \gamma_0(t)$ lies in the kernel of  the restriction of $\omega(\lambda)$ on the fiber $F(\gamma_0(t))$:
\begin{equation} \label{ker2}
\omega(\lambda)(W_k(\gamma_0(t)), \dot \gamma_0(t))=0,\ \ k=1,\dots,n.
\end{equation}
It is straightforward to compute that $\omega_i(W_k(\gamma_0(t)), \dot \gamma_0(t))=0$ unless $i=k$ or $i=k+1$; in these two cases, one has:
$$
\omega_k(W_k(\gamma_0(t)), \dot \gamma_0(t))=C,\ \ \omega_{k+1}(W_k(\gamma_0(t)), \dot \gamma_0(t))=-C
$$
where $C$ is a non-zero constant depending on $n$ only. Thus (\ref{ker1}) and (\ref{ker2}) imply that $2\lambda_k-\lambda_{k-1}-\lambda_{k+1}=0$ for all $k$. This, in turn, implies that all $\lambda_k$ are equal (indeed, by induction, $\lambda_{k+1}=k\lambda_2 - (k-1)\lambda_1$ for all $k$, which implies, for $k=n$, that $\lambda_1=\lambda_2$). Since $\sum \lambda_i=0$, all  $\lambda_k$ vanish. It follows that  the curve $\gamma_0$ is not singular.

Finally, if the rotation number $k$ of a periodic trajectory is greater than 1, one changes the curve (\ref{1curve}) to
$$
Z(t)=(z_1(t),\dots,z_n(t)) \ \ {\rm with}\  \ z_j(t)=e^{2\pi i(jk+t)/n},\ t\in[0,1],
$$
and the proof proceeds along the same lines. 
\proofend

\begin{remark} \label{sing}
{\rm To illustrate the importance of checking that the curve $\gamma_0$ is not singular, let  us give a simple example of a rigid horizontal curve, that is, a curve that does not admit smooth perturbations in the class of horizontal curves with fixed end-points, see \cite{Mon}.  Consider  the curve $\gamma(t)=(t,0,0),\ t\in[0,1]$ tangent to the distribution $dz=y^2dx$ in 3-space. If   $(x(t),y(t),z(t))$ is its perturbation as a horizontal curve with fixed end-points then
$$
0=z(1)-z(0)=\int_0^1 y^2(t) x'(t) dt.
$$
Since $x'(t)>0$, we have $y(t)=0$, and hence $z(t)=0$, for all $t$; thus the perturbed curve is a reparameterization of $\gamma$.
}
\end{remark}

\section{Case study:  $n=3$ and discretization} \label{per3}

Three is the smallest possible period of the outer billiard map. In this section we 
illustrate the general constructions of Sections \ref{conf} and  \ref{perpts} in the case $n=3$.

Let $\gamma$ be a smooth strictly convex outer billiard curve and $z_1 z_2 z_3$ a 3-periodic orbit. Let $w_1, w_2, w_3$ be the mid-points of the sides of the triangle $z_1 z_2 z_3$. 

\begin{lemma} \label{inout}
The  inscribed triangle $w_1 w_2 w_3$ has the property that, for every  $i=1,2,3$, the tangent line to $\gamma$ at vertex $w_i$ is parallel to the side $w_{i+1} w_{i+2}$. Conversely, if an inscribed triangle has this property then the tangent lines to $\gamma$ at its vertices form a circumscribed triangle whose sides are bisected by the tangency points. Furthermore, this property of an inscribed triangle is equivalent to having extremal area.
\end{lemma}

\proof The first claim holds by definition of the outer billiard map and the second follows from elementary geometry. To prove the last statement, fix points $w_2$ and $w_3$ and vary $w_1$. The area is extremal when  the distance from $w_1$ to the line $w_2 w_3$ has an extremum, and this happens if and only if the tangent line to $\gamma$ at $w_1$ is parallel to  $w_2 w_3$.
\proofend

It follows that we can replace circumscribed triangles of extremal areas by inscribed ones. We want to construct smooth strictly convex closed curves such that an inscribed triangle of extremal area can be turned around inside the curve.

Following the approach of Section \ref{conf}, let $M$ be the 5-dimensional manifold of triangles of area 3/2. Let $w_1 w_2 w_3$ be a triangle from $M$; denote by $c$ its center of mass. Let $u=w_1-c, v=w_2-c$; then $[u,v]=1$ and the frame $(u,v)$ belongs to $SL(2,\R)$.  The cyclic group $\Z_3$ acts on $M$ by the cyclic permutations of the vertices of a triangle and on $SL(2,\R)$ by the transformation
$T(u,v)=(v,-u-v). $
Extend the action of $\Z_3$ to $SL(2,\R)\times \R^2$ taking product with the trivial action on the second factor. On $M$, we have a 3-dimensional distribution ${\cal F}$ determined by the condition that the velocity of $w_i$ is parallel to the line $w_{i+1} w_{i+2},\ i=1,2,3$.

\begin{lemma} \label{decomp}
The correspondence $(w_1, w_2, w_3) \mapsto ((u,v),c)$ is a $\Z_3$-equivariant diffeomorphism  $M\to SL(2,\R)\times \R^2$. The fibers of the projection $SL(2,\R)\times \R^2 \to SL(2,\R)$ are transverse to the distribution ${\cal F}$.
\end{lemma}

\proof The first claim is obvious. To prove the second, note that a vector, tangent to the fibers of the projection $SL(2,\R)\times \R^2 \to SL(2,\R)$, is induced by a parallel translation of the triangle: all three vertices move with the same velocity. If this vector  lies in ${\cal F}$ then it must be parallel to all three sides of the triangle, which is impossible.
\proofend

We want to construct a closed horizontal curve $\tilde \gamma(t),\ t\in[0,2\pi]$ in $M$ satisfying the monodromy condition 
\begin{equation} \label{monodr}
\tilde\gamma(t+2\pi/3)=\sigma(\tilde\gamma(t))
\end{equation}
 where $\sigma$ is the cyclic permutation of the vertices of a triangle. By Lemma \ref{decomp}, this curve projects to a smooth curve $\gamma(t)$ in $SL(2,\R)$ satisfying  $\gamma(t+2\pi/3)=T(\gamma(t))$. Conversely, such a curve $\gamma(t)\subset SL(2,\R)$ lifts uniquely to a horizontal curve $\tilde \gamma(t) \subset M$, provided a lift of the initial point is chosen, but the monodromy condition (\ref{monodr}) may fail. But once (\ref{monodr}) is satisfied -- and assuming convexity -- we obtain a desired plane curve possessing a 1-parameter family of inscribed triangles of extremal areas.
 
Let us work out explicit formulas. Denote by $\gamma(t)=(u(t),v(t)),\ t\in[0,2\pi]$ a closed curve in $SL(2,\R)$ such that $u(t+2\pi/3)=v(t)$ and $v(t+2\pi/3)=-u(t)-v(t)$. Let $c(t)$ be the curve of the center of mass of the respective triangle, that is, the projection to $\R^2$ of a lifted horizontal curve $\tilde \gamma(t)$. Set:
$$
[u(t),u'(t)]=p(t),\ [v(t),v'(t)]=q(t),\ [u'(t),v(t)]=[v'(t),u(t)]=r(t).
$$  
 
\begin{proposition} \label{form}
The velocity of the center of mass is given by the formula: 
\begin{equation} \label{cm}
c'(t)=\frac{1}{3}(p(t)-2q(t)+2r(t))\ v(t)-\frac{1}{3}(q(t)-2p(t)+2r(t))\ u(t).
\end{equation}
 The monodromy condition (\ref{monodr}) holds if and only if 
\begin{equation} \label{mon2}
\int_0^{2\pi} (p(t)-2q(t)+2r(t))\ v(t)-(q(t)-2p(t)+2r(t))\ u(t)\ dt =0
\end{equation}
or, equivalently,
\begin{equation} \label{mon3}
\int_0^{2\pi} [u'(t),v(t)] (v(t)-u(t)) \ dt=0.
\end{equation}
\end{proposition}

\proof One has:
\begin{equation} \label{wc}
w_1(t)=u(t)+c(t),\ w_2(t)=v(t)+c(t),\ w_3(t)=-u(t)-v(t)+c(t).
\end{equation}
The curve is horizontal if and only if
$$
[w_2(t)-w_1(t),w'_3(t)]=[w_3(t)-w_2(t),w'_1(t)]=[w_1(t)-w_3(t),w'_2(t)]=0
$$
or, in view of (\ref{wc}),
$$
[u(t),c'(t)]+2[v(t),c'(t)]=2r(t)-p(t), 2[u(t),c'(t)]+[v(t),c'(t)]=2r(t)-q(t)
$$
(there are only two equations because the third is their consequence). It follows that
$$
3[u(t),c'(t)]=2r(t)-2p(t)+q(t),\ 3[v(t),c'(t)]=2r(t)-2q(t)+p(t),
$$
and since $[u(t),v(t)]=1$, this implies (\ref{cm}).

Abusing notation, let $T$ be the shift $t\mapsto t+2\pi/3$ of the argument of the functions $u(t),v(t),p(t),q(t),r(t)$; this defines an action of $\Z_3$ on the space generated by these functions. One has:
\begin{equation} \label{act}
T:u\mapsto v\mapsto -u-v,\ p\mapsto q\mapsto p+q-2r,\ r\mapsto -r+q \mapsto p-r.
\end{equation}
The monodromy condition (\ref{monodr}) is satisfied if and only if $c(t)$ is $2\pi/3$-periodic. Applying (\ref{act}) to  (\ref{cm}), it follows that $c'(t)$ is $2\pi/3$-periodic, and hence  $c(t)$ is $2\pi/3$-periodic if and only if $c(t)$ is $2\pi$-periodic, that is, if $\int_0^{2\pi} c'(t) dt =0$. In view of (\ref{cm}), this is equation (\ref{mon2}). Finally, for every function $f(t)$, one has
$$
\int_0^{2\pi} T(f)(t)\ dt = \int_0^{2\pi} f(t)\ dt.
$$
Applying this to (\ref{act}) yields four identities for integrals of vector-functions $pu,pv,qu,qv,ru,rv$:
$$
\int_0^{2\pi} q(t)u(t)\ dt =-\int_0^{2\pi} r(t) v(t) \ dt, \int_0^{2\pi} p(t) v(t)\ dt=-\int_0^{2\pi} r(t) u(t)\ dt,
$$ 
$$
\int_0^{2\pi} q(t) v(t)\ dt=\int_0^{2\pi} p(t) u(t)\ dt=\int_0^{2\pi} r(t) u(t)\ dt+\int_0^{2\pi} r(t) v(t)\ dt.
$$
 These identities imply that (\ref{mon2}) and (\ref{mon3}) are equivalent.
\proofend

To summarize, one constructs an outer billiard table having a 1-parameter family of 3-periodic trajectories from a closed curve in $SL(2,\R)$ satisfying  (\ref{mon3}), that is, satisfying two scalar equations. 

It is interesting to consider the following discretization of the problem in which a smooth curve is replaced by a polygon. Let $P$ be a convex polygon and $ABC$ an inscribed triangle whose vertices are some vertices of $P$. One is allowed to slide one vertex of the triangle along an edge of $P$ to the adjacent vertex of $P$,  provided this edge is parallel to the opposite side of the triangle. Using a sequence of such moves, one wants to turn the triangle around inside $P$. For which polygons is this possible?

\begin{figure}[hbtp]
\centering
\includegraphics[width=5in]{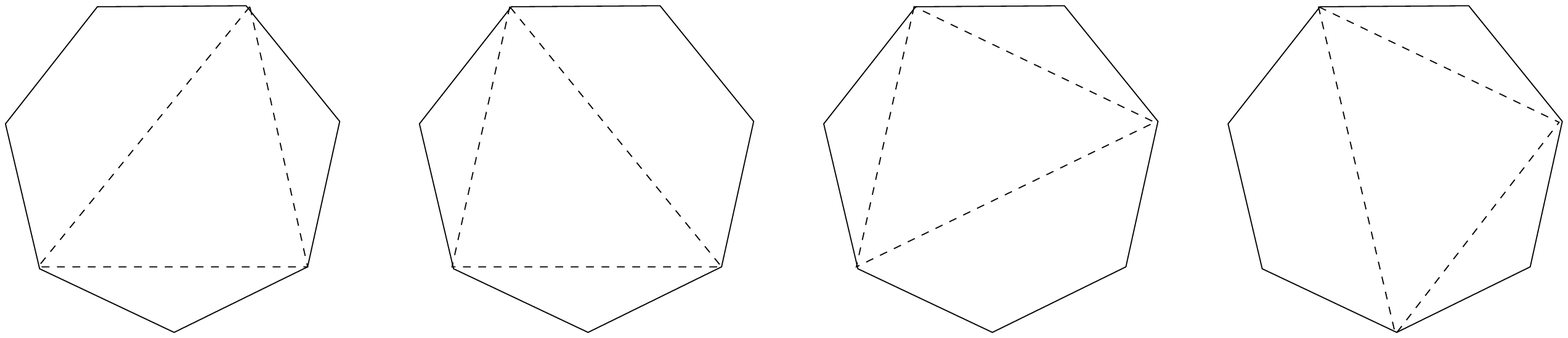}
\caption{Rotating a triangle inside a regular septagon}
\label{septa}
\end{figure}

Consider an example. Let $P$ be a regular $k=3n+1$-gon with vertices $z_1,\dots z_{k}$, and let $A=z_1, B=z_{n+2}, C=z_{2n+2}$. Then the segment $z_1 z_2$ is parallel to the segment $BC$, and one can move point $A$ to $z_2$. After that one can move point  $C$ to  $z_{2n+3}$, then point $B$ to $z_{n+3}$, etc. -- see figure \ref{septa} for $k=7$.

One can deform a regular $k=3n+1$-gon preserving the desired property: the condition to satisfy is that the segments $z_i z_{i+1}$ and $z_{i+n+1} z_{i+2n+1}$ be parallel for all $i=1,\dots, k$. This gives $k$ conditions, whereas a $k$-gon has $2k$ degrees of freedom, so it is not surprising that non-trivial deformations exist. We do not dwell on details here.

Let us finish this section with a problem. {\it Suppose that a convex domain has the property that an inscribed triangle of extremal area can be turned around inside the domain. Find the lower and upper bounds for the ratio of the areas of the triangle and the domain.} We conjecture that a parallelogram and an ellipse are two extreme cases.

\section{The set of 3-periodic orbits has empty interior} \label{zerom}

It is an old conjecture, with motivation in spectral geometry, that the set of $n$-periodic billiard trajectories in a billiard  with a smooth boundary (not necessarily convex) has zero measure. For $n=2$, this is easy to show, and for $n=3$, this is a theorem of Rychlik \cite{Ry} with a number of different proofs available \cite{BZh,Sto,Vo,Wo}. For $n\geq 4$, the conjecture is open.

This section concerns the outer counterpart of this conjecture. Consider $n$ outer billiard tables $D_1,\dots,D_n$  with smooth boundaries and let $T_1,\dots,T_n$ be the respective outer billiard maps. We are interested in the  set $F$ of fixed points of the composition $T_n\circ T_{n-1}\circ \dots \circ T_1$.

\begin{conjecture} \label{conjzero}
For every $n\geq 3$, the set $F$ has zero measure.
\end{conjecture}

\begin{figure}[hbtp]
\centering
\includegraphics[width=2in]{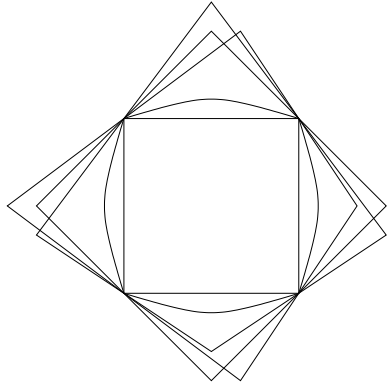}
\caption{An open set of 4-periodic orbits for a piece-wise smooth outer billiard table}
\label{round}
\end{figure}

\begin{remark} \label{conjrem}
{\rm  For strictly convex outer billiards with corners, Conjecture \ref{conjzero} fails. An example shown in figure \ref{round} is a rounded square: the set of 4-periodic orbits has a non-empty interior.
One can generalize: if a plane polygonal outer billiard has an $n$-periodic orbit  then a whole neighborhood of this point consist of periodic points as well (with period $2n$, if $n$ is odd). According to a theorem of  Culter \cite{Cu}, every polygonal outer billiard has periodic orbits,\footnote{For inner polygonal billiards, this is an outstanding open conjecture.} thus one can use any (rounded) polygon to construct an open set of periodic orbits.
}
\end{remark} 

We prove the following particular case  of Conjecture \ref{conjzero}.

\begin{theorem} \label{3per}
Let $D_1,D_2,D_3$ be three smooth strictly convex outer billiard tables and $T_1,T_2,T_3$ the respective outer billiard maps. Then the set $F$ of fixed points of the map $T_3\circ T_2\circ T_1$ has empty interior.
\end{theorem}

As in \cite{Ry,Sto,Vo,Wo}, this theorem implies that  $F$ has zero measure; we do not dwell on this refinement.

As a preparation to the proof of Theorem \ref{3per}, let us recall two results on outer billiards. The first is an outer analog of the ``mirror equation" of geometrical optics that plays an important role in the theory of inner billiards. 

\begin{figure}[hbtp]
\centering
\includegraphics[width=3in]{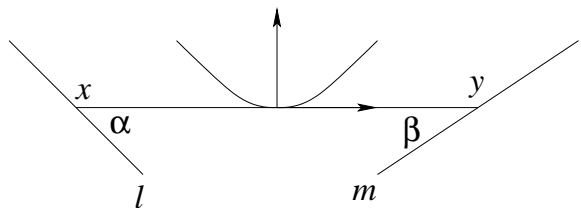}
\caption{The differential of the outer billiard map and the mirror equation}
\label{dif}
\end{figure}

Let $D$ be an outer billiard table with smooth boundary and $T$ the outer billiard map.
Let $T(x)=y$, see figure \ref{dif}.

\begin{proposition} \label{mirror}
In the Cartesian coordinate system in figure \ref{dif}, the differential of $T$ at point $x$ is given  by the matrix
$$
 dT=\left(\begin{array}{cc}
 -1&-\frac{2\rho}{r}\\
 0&-1
 \end{array}\right)
$$
where $r=|xy|/2$ and $\rho$ is the radius of curvature of the outer billiard curve at the reflection point. If $dT$ takes a line $l$ to a line $m$, see figure \ref{dif}, then
one has the following  ``mirror equation":
\begin{equation} \label{mireq}
\cot \alpha + \cot \beta = \frac{2\rho}{r}.
\end{equation}
\end{proposition}

For a proof of Proposition \ref{mirror}, see \cite{G-K}. 

The second result is the following  ``area construction", an outer counterpart of the string construction in the theory of inner billiards, that recovers a billiard table from its caustic.

\begin{proposition} \label{area}
Let $D$ be an outer billiard table and $C$ an invariant curve of the outer billiard map $T$. Let $x\in C$ and $T(x)=y$. Then the area cut by the segment $xy$ from $C$ is the same for all $x\in C$. In other words, $D$ can be recovered from $C$ as an envelope of segments of constant areas.
\end{proposition}

For a proof of Proposition \ref{area}, see \cite{D-T,Tab1,Tab2} or \cite{G-K}.

In particular, let $C$ consist of two non-parallel lines, $l$ and $m$, as in figure \ref{dif}. Then the area construction  yields a hyperbola for which the lines $l$ and $m$ are the asymptotes. If the outer billiard map $T$ takes $x$ to $y$ and  $dT(l)=m$, as in figure \ref{dif}, then the hyperbola  for which the lines $l$ and $m$ are the asymptotes is second order tangent to the outer billiard curve at the mid-point of the segment $xy$.
Note also that of 4 possible orientations of the pair of lines $l$ and $m$ in figure \ref{dif}, only two are consistent with the action of $dT$: one line is oriented toward their intersection point while the other is oriented from it.

\paragraph{Proof of Theorem \ref{3per}.}
We give two proofs. The first is reminiscent of the one given in \cite{Wo}. 

\begin{figure}[hbtp]
\centering
\includegraphics[width=3in]{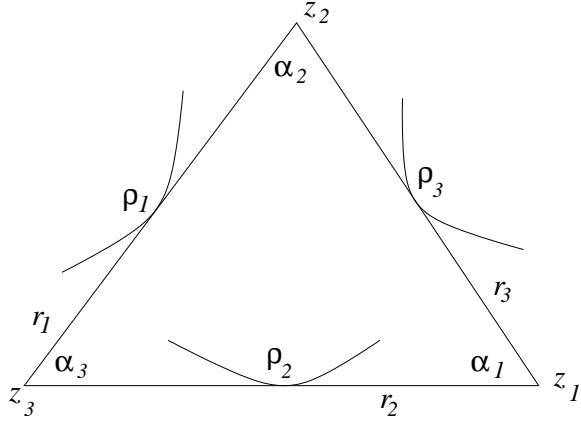}
\caption{Proving Theorem \ref{3per}}
\label{3pf}
\end{figure}

Consider figure \ref{3pf}. One has three outer billiard maps, $T_1,T_2,T_3$ that cyclically permute points $z_1,z_2,z_3$:
$$
T_1(z_2)=z_3,\ T_2(z_3)=z_1,\ T_3(z_1)=z_2.
$$
 Let $\alpha_i,\ i=1,2,3$, be the angles of the triangle, $2r_i$ its side lengths and $\rho_i$ the curvature radii of the outer billiard curves at the reflection points. Set $dT_i=A_i$ and let $B_i$ be matrix of rotation through angle $\pi-\alpha_i$.

Assume that the orbit $z_1,z_2,z_3$ belongs to an open set $U$ of fixed points of the map $T_3\circ T_2\circ T_1$. Then $d(T_3\circ T_2\circ T_1)$ is the identity in $U$. Using Proposition \ref{mirror}, this implies that 
\begin{equation} \label{mat}
B_1A_2B_3A_1B_2A_3={\rm Id}\quad {\rm or}\quad A_1B_2A_3=B_3^{-1}A_2^{-1}B_1^{-1}.
\end{equation}
The matrices involved are as follows:
$$
A_i=\left(\begin{array}{cc}
 -1&\frac{-2\rho_i}{r_i}\\
 0&-1
 \end{array}\right),
 \quad
 B_i=\left(\begin{array}{cc}
 -\cos \alpha_i&-\sin \alpha_i\\
 \sin \alpha_i&-\cos \alpha_i
 \end{array}\right).
$$
It is straightforward to compute the triple products of the matrices on both sides of the second equation (\ref{mat}). Apply both products to the vector $(1,0)$ and equate the second components:
$$
\sin \alpha_2 = -\sin (\alpha_1+\alpha_3) + 2\sin \alpha_1 \sin \alpha_3 \frac{\rho_2}{r_2},
$$
and since $\alpha_2=\pi-\alpha_1-\alpha_3$,
\begin{equation} \label{expr}
\frac{\rho_2}{r_2}= \frac{\sin (\alpha_1+\alpha_3)}{\sin \alpha_1 \sin \alpha_3}=\cot \alpha_1 + \cot \alpha_3.
\end{equation}

Let $A$ be the area of the triangle $z_1,z_2,z_3$ and $h_2$ its altitude from the vertex $z_2$. Then by elementary geometry,
$$
2r_1 \sin \alpha_3=2r_3 \sin \alpha_1=h_2,\ \ r_1 \cos \alpha_3+r_3 \cos \alpha_1=r_2,\ \ A=h_2r_2,
$$ 
which makes it possible to rewrite (\ref{expr}) as
\begin{equation} \label{expr1}
A \rho_2=2r_2^3.
\end{equation}

For  orbits $z_1,z_2,z_3$ in the set $U$, the area $A$ remains constant. Move points 
$z_1$ and $z_3$ slightly along the line $z_1z_3$ so that the mid-point of the segment $z_1 z_3$ remain the same. We obtain a family of orbits in $U$ for which $\rho_2$ does not change but $r_2$ does, in contradiction with (\ref{expr1}). This completes our first proof.

The second argument is a variation on the proof given in \cite{BZh}. Let $M^5$ be the manifold of non-degenerate triangles whose oriented areas equal 1 and let  ${\cal F}^3$ be the dual Birkhoff distribution on $M$. If there is an open set of fixed points of the map $T_3\circ T_2\circ T_1$ then one has a 2-dimensional disc  $U\subset M$ tangent to ${\cal F}$. We will show that this is impossible.

If such a disc $U$ exists, let  $x$ and $y$ be local coordinates in it. Then $u=\partial/\partial x$ and $v=\partial/\partial y$ are commuting  vector fields on $U$.

As a basis of vector fields in $M$ generating ${\cal F}$ we can choose the fields
$$
W_1=(z_2-z_1)\partial z_2 + (z_1-z_3)\partial z_3,\ \  W_2=(z_3-z_2)\partial z_3 + (z_2-z_1)\partial z_1
$$
and
$$
 W_3=(z_1-z_3)\partial z_1 + (z_3-z_2)\partial z_2,
$$
see Section \ref{conf}.  Then 
$$
u=f_1 W_1+f_2 W_2+f_3W_3,\ \ v=g_1 W_1+g_2 W_2+g_3W_3
$$
where $f_i,g_i$ are smooth functions on $U$. Without loss of generality, $f_1\neq 0$ and $g_2\neq 0$. Then we take a linear combination of $u$ and $v$, with coefficients that are smooth functions, to obtain new fields
$$
\bar u=W_1+fW_3,\ \ \bar v=W_2+gW_3,
$$
so that $[\bar u,\bar v] \in {\rm Span}(\bar u,\bar v)$. In particular, $[\bar u,\bar v]$ is tangent to ${\cal F}$.

It follows that $d\alpha_i (\bar u, \bar v)=0$ for $i=1,2,3$ where  $\alpha_i$ are the 1-forms introduced in Section \ref{conf}. It is straightforward to compute  that
$$
d\alpha_2(\bar u, \bar v)=-2(1+g),\ \ d\alpha_3(\bar u, \bar v)=2(1+f),
$$
and hence $f=g=-1$. Thus $\bar u=W_1-W_3, \bar v=W_2-W_3$.

Now one computes the commutator:
$$
[W_1-W_3,W_2-W_3]=[W_1,W_2]+[W_2,W_3]+[W_3,W_1]=4(W_1+W_2+W_3).
$$
Since the fields $W_1,W_2,W_3$ are everywhere linearly independent, the vector $W_1+W_2+W_3$ is not a linear combination of $W_1-W_3$ and $W_2-W_3$. This is a contradiction proving that the disc $U$ does not exist.
\proofend

\begin{remark} \label{difl}
{\rm The reason we consider $n$ outer billiards in Conjecture \ref{conjzero}, instead of just one, is that this makes it possible to relax the convexity assumptions and include outer billiard tables, convex outward (similarly to the case of inner billiards where one does not assume the billiard table to be convex). 

In fact, the proof of Theorem \ref{3per} considerably simplifies if one considers a single strictly convex outer billiard table. Namely,  let $z_1$ be a 3-periodic point of an outer billiard map $T$. Then the differential $dT^3 (z_1)$ cannot be equal to identity. 

\begin{figure}[hbtp]
\centering
\includegraphics[width=5in]{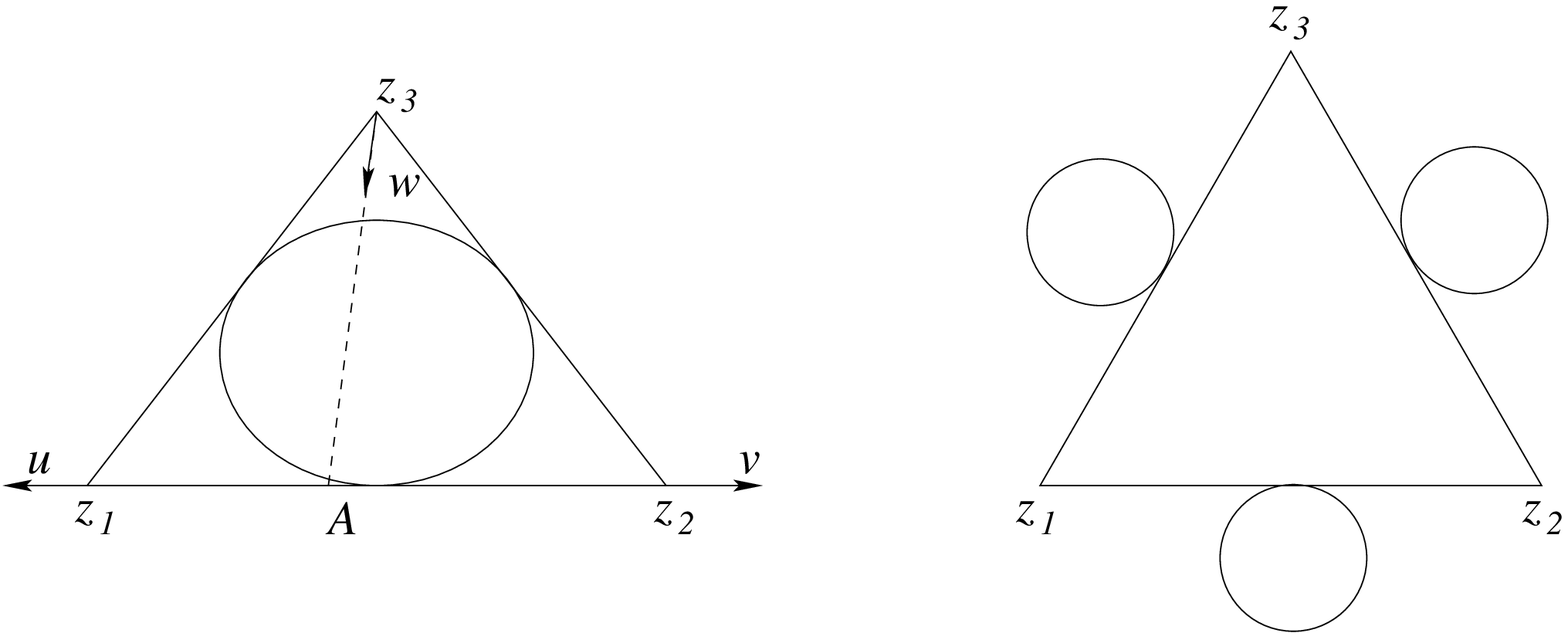}
\caption{Can $dT^3={\rm Id}$?}
\label{exa}
\end{figure}

Indeed,  let $T(z_1)=z_2, T(z_2)=z_3, T(z_3)=z_1$, see figure \ref{exa}. Assume that $dT^3 (z_1)={\rm Id}$. Let $u$ and $v$ be the opposite unit tangent vectors at points $z_1$ and $z_2$ having the direction of the side $z_1z_2$. Then $dT(u)=v$. Let $w=dT(v)=dT^{-1}(u)$. Since $T$ is an orientation preserving transformation, the vector $w$ lies on the same sides of the lines $z_2z_3$ and  $z_1z_3$ as the outer billiard table $D$. Thus $w$ lies inside the wedge $z_1z_3z_2$. Let $A$ be the intersection point of the line through point $z_3$ in the direction of vector $w$ with the line  $z_1z_2$. It follows from the area construction (see the paragraph after Proposition \ref{area}) that the outer billiard curve is second order tangent to the hyperbola inscribed in the angle $z_1Az_3$ and tangent to the segment $z_1z_3$ in its mid-point. This hyperbola is convex outward, contradicting the convexity of $D$.

On the other hand,  $d(T_2\circ T_1\circ T_3)$ may be equal to the identity at a fixed point of the map $T_2\circ T_1\circ T_3$. An example is shown in figure \ref{exa} on the right: $z_1z_2z_3$ is a unit equilateral triangle, and the three outer billiard tables are circles of radii $1/\sqrt{3}$. Then equation (\ref{expr1}) holds, and $d(T_2\circ T_1\circ T_3)(z_1)={\rm Id}$.
}
\end{remark}

\medskip

{\bf Acknowledgments.} We are grateful to Yu. Baryshnikov, S. Lvovski, R. Montgomery, K. F. Siburg, V. Zharnitskii for stimulating discussions. The second author is indebted to Max-Planck-Institut in Bonn for its invariable hospitality.

\end{document}